\documentclass[11pt]{article}

\usepackage{amsmath,amssymb,amsthm,amscd,graphics,psfrag,fancyheadings,float}
\usepackage{pstricks, pst-plot, epsfig}
\usepackage{stmaryrd}
\usepackage{caption}
\usepackage{refstyle}

\theoremstyle{plain}

\newtheorem{theorem}{Theorem}[section]

\theoremstyle{definition}

\newcommand{\C}{{\mathbb C}} 
\newcommand{\N}{{\mathbb N}}

\newcommand{\cale}{{\mathcal E}}

\divide\oddsidemargin by 2
        {\end{list}} 

\newsavebox{\savepar}

\begin{document}

\title{{\Large A  complementary proof   of  Baker's  theorem  of completely invariant components for transcendental  entire functions}}

\author{
P. Dom\'{\i}nguez$^a$ and  G. Sienra$^b$ \\ {\small  (a) Fac. Ciencias F\'{\i}sico-Matem\'aticas, Benem\'erita Universidad Aut\'onoma de Puebla} \\ {\small   Avenida San Claudio y 18 Sur, C.U., Puebla Pue, 72570, M\'exico.}\\
 {\small  (b) Facultad de Ciencias, Universidad  Nacional Aut\'onoma de M\'exico}\\  {\small Avenida Universidad 3000, C.U. Ciudad de M\'exico, 04510, M\'exico}}
\date{}

\maketitle
      
\begin{abstract}   Baker in \cite{baker1e}  proved that  for transcendental entire functions  there is 
at most one completely invariant component of the Fatou set.   It was  observed by Julien Duval  that there is a missing case  in Baker's proof. In this article we follow Baker's ideas and  give some  alternative arguments to establish the result.

\renewcommand{\thefootnote}{} 
 \footnote{2017  
{\it pdsoto@fcfm.buap.mx , guillermo.sienra@gmail.com}} 
\addtocounter{footnote}{-2}            
\end{abstract}
 
\section{Introduction}

Let $\cale$ be  the set of  trascendental entire  functions $f : \mathbb{C} \rightarrow \mathbb{C}$. For $f \in \cale$, we write $f^n=f\circ f^{n-1}$ for the \textit{n}-th iterate of $f$, $n \in \mathbb{N}$, and $f^0= Id$ where the symbol $\circ$ denotes composition. When $f^n(z_0)=z_0$, for some $n\in \mathbb{N}$,  the point $z_0$ is called a periodic point  If $n$ is the minimal positive integer for which this equality holds, we say that $z_0$ has period $n$. If $n=1$, $z_0$ is called a fixed point.  The classification of a periodic point $z_0$ of period $n$ of $f \in \cale$ can be  attracting, super-attracting, rationally indifferent, irrationally indifferent and repelling.

Given $f \in \cale$, the \textit{Fatou set} $\mathfrak{F}(f)$ is defined as the set of all points $z \in\mathbb{C}$ such that the sequence of iterates $(f^n)_{n\in \mathbb{N}}$ forms a normal family in some neighborhood of $z$. The \textit{Julia set}, denoted by $\mathfrak{J}(f)$, is the complement of the Fatou set. 

Some properties of the Julia and Fatou sets for functions in class $\cale$  are mentioned below: 

(i) $\mathfrak{F}(f)$ is open, so $\mathfrak{J}(f)$ is  closed.

(ii)  $\mathfrak{J}(f)$ is perfect and non-empty.

 (iii) The sets $\mathfrak{J}(f)$ and $\mathfrak{F}(f)$ are completely  invariant under $f$. 
 
  (iv) $\mathfrak{F}(f)=\mathfrak{F}(f^n)$ and $\mathfrak{J}(f)=\mathfrak{J}(f^n)$ for all $n \in \mathbb{N}$.
  
   (v) The repelling periodic points are dense in $\mathfrak{J}(f)$. \\
   
   See \cite{berg}, \cite{eremenko1},  \cite{eremenko2} and \cite{eremenko3} for definitions,  proofs and more details concerning the Fatou and Julia sets. \\

We denote by $CV$ the set of critical values and  by $OV$ the set of  omitted values of a function $f \in \cale$.  

We recall that a Fatou component $G$ of $f$ is {\sl completely invariant} if $f^{-1}(G)=G$. Also, for any two points $z_{1},z_{2} \in G$, there is a path contained in $G$ that joins the two points. For $f$ a transcendental entire function, due to Picard theorem, every completely invariant Fatou component of $f$ is unbounded. If $f$ has $k$ completely invariant components, $G_{k}$, with $k \in \N$, then for every point $z \in \mathfrak{J}(f)$ and any neighborhood $N_{z}$ of $z$, we have $G_{k} \bigcap N_{z} \neq \emptyset$. \\

\noindent {\bf Observation 1.}  {\it Let $f \in \cale$, and $G$ a completely invariant Fatou component of $f$. Let $w \in G$ a regular value of $f$ and $z(w)$ any pre-image point, then there exist an oriented curve $\Gamma \subset G$ beginning at $w$,  such that: (i) intersects any neighborhood of infinity with  $\Gamma \bigcap OV=\emptyset$ and (ii) has a pre-image $\Gamma'$ beginning at $z(w)$, so $f(\Gamma')=\Gamma$.

The  construction of the curve $\Gamma$ can be obtained by successive applications of a generalization of the Gross-star theorem, due to Kaplan \cite{ka} Theorem 3. In few words, Kaplan proves in particular, that for a (star) family of non intersecting bounded curves begining at a regular value $w$, the pre-images based at any $z(w)$ exist and can be continued indefinitely, for almost all of the curves, see details in  \cite{ka} (compare with Inversen's Theorem  in \cite{eremenko4}).

To construct $\Gamma$ we   consider any neighborhood $N_{1}$ of $\infty$  and any pre-image $w_{1}$ of $w$ in $N_{1}$, then let $\tau_{1}$ any path contained in $G$ joining $w$ and $w_{1}$ (with that orientation). By Kaplan's theorem, for any small enough neighborhood $B(w_{1})$ of $w_{1}$, there is a curve $\tau_{1}'$ beginning at $w$ and ending in some $w_{1}' \subset B(w_{1})$ which has a pre-image $\tau_{2}$, beginning at $w_{1}'$. 

Now, proceed inductively by choosing  neighborhoods $N_{i} \subset N_{i-1}$, with $N_{i}$ converging to $\infty$ as $i$ tends to $\infty$ and choosing points $w_{i} \subset N_{i}$, with $w_{i}$ a pre-image of $w_{i-1}'$, also take paths $\tau_{i}$   joining $w_{i-1}'$ with $w_{i}$ and modify them to $\tau_{i}'$ accordingly to Kaplan's theorem. \\

In conclusion $\Gamma:=\bigcup_{i=1}^{\infty} \tau_{i}'$  is a curve satisfying  (i) and (ii)  and  its closure $\overline{\Gamma}$ is a continuum  in the sphere.}

An important result related to completely invariant components   of transcendental entire functions was given by Baker in \cite{baker1e}, it is stated as follows.

\begin{theorem}
If $f \in \cale$, then there is 
at most  one completely invariant component of $\mathfrak{F}(f)$.
\label{ba}
\end{theorem}

As it was mentioned in the abstract   that there is a missing case in Baker's proof,  in this paper  we follow Baker's ideas  and  give some alternative arguments  to solve the missing case.\\

It is interesting to note that a  recent paper by  Rempe and Sixmith \cite{rpe},  studies the connectivity of the pre-images of simply-connected domains of a  transcendental  entire function. The  paper  describes in detail the error in Baker's proof and mentions  Duval's  example, which is equivalent to the case of Figure \ref{fig4} in this article. Also,  they prove that if infinity is accesible from some Fatou component, then at least one of the pre-images of some component is disconnected. Since infinity is accesible in Baker domains, they conclude that if the function  has two completely invariant Fatou components both components must be attracting or parabolic basins. In their article  it is included a list of papers which use  Baker's result. 

While we  were making  final corrections of  this paper we got,   by communication, the results obtained by  Rempe and Sixmith in \cite{rpe}.

\section{Proof of Theorem} 

The idea of the proof is by contradiction assuming  that there are at least two completely invariant open components $G_{1}$ and  $G_{2}$. We begin  considering the cut system of Baker  in Step 1, which is an open disc $D_{1}$ with boundary the simple curves $\widehat\gamma_{1}$, $\beta_{1}, \widehat\gamma_{2}, \beta_{2}$, with the properties that $\beta_{1} \subset G_{1}$ and $\beta_{2} \subset G_{2}$ and such that $f(\widehat\gamma_{1}) \subset \gamma$ , $f(\widehat\gamma_{2}) \subset \gamma$ are conformal injections, for $\gamma$ a segment with extremes at $G_{1}$ and $G_{2}$. In Step 2, we consider the image of the disc $f(D_{1})$ which has to be bounded and state some of its properties. We proceed in Step 3 to extend the curves $f(\beta_{1})$ and $f(\beta_{2})$ to infinity as in the Observation 1, creating two unbounded curves $\Gamma \subset G_{1}$ and $\Theta \subset G_{2}$. Such curves can be very complicated inside $f(D_{1})$, so we consider their intersection with the complement of $f(D_{1})$ that we named $B$. Then, we studied their pre-images in the complement of $D_{1}$. By adding a certain path $\sigma$ (Case A) and $\Sigma$ (Case B) between those pre-images in the same component, we show that one of the regions $G_{1}$ or $G_{2}$ is disconnected, which is a contradiction. It is important for the proof, the cut system since it helps to have certain control on the pre-images of  $\Gamma \bigcap B = \Gamma_0$ and $\Theta \bigcap B = \Theta_0$. The differences between Case A and Case B relay in the way the pre-images of $\Gamma_0$ and $\Theta_0$ intersects the cut system, as indicated in the Step 3.\\

\noindent {\bf Proof.} Suppose that $\mathfrak{F}(f)$ has at least two mutually disjoint  completely invariant 
domains  $G_1$ and $G_2$.\\

\noindent {\bf Step 1. The Cut System and the cancelation procedure.}

\begin{figure}[h!]
\centerline{\hbox{\psfig{figure=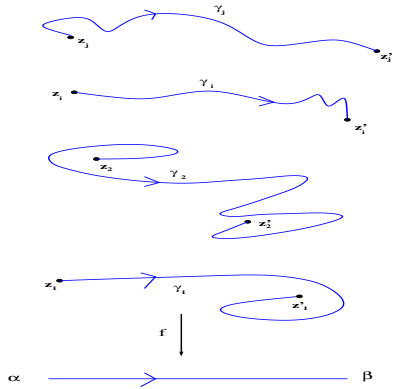,height=2in,width=2in}}}
\caption{The  images $g_i(\gamma) = \gamma_i$}
\label{fig0}
\end{figure}

Take a value $\alpha$ in $G_1$ such that $f(z) = \alpha$ has infinitely many 
simple roots $z_i$ \, ($f'(z) = 0$ at only countably many $z$ so we have to 
avoid only countably many choices of $\alpha$). All $z_i$ are in $G_1$. Similarly take $\beta$ in $G_2$  such that $f(z) = \beta$ has infinitely many simple roots $z_i'$ in $G_2$. By Gross' star theorem 
\cite{gross} we can continue 
all the regular branches $g_i$ of $f^{-1}$ such that $g_i(\alpha) = z_i$, 
along almost every ray to $\infty$ without meeting any singularity (even algebraic).
 Thus we can move $\beta \in G_2$ slightly if necessary so that all $g_i$ 
continue to $\beta$ analytically along the line $\gamma$, which joins $\alpha$ 
and $\beta$. The images $g_i(\gamma)$ are disjoint curves joining $z_i$ to $z_i'$.
Denote $g_i(\gamma) = \gamma_i$. Note that $\gamma_i$ is oriented from 
$z_i$ to $z_i'$, see Figure \ref{fig0}.

The branches $f^{-1}$ are univalent so $\gamma_i$ are disjoint simple arcs. 
Different $\gamma_i$ are disjoint since $\gamma_i$ meets $\gamma_j$ at
 say $w_0$ only if two different 
branches of $f^{-1}$ become equal with values $w_0$ which can occur only if 
$f^{-1}$ has branch point at $f(w_0)$ in $\gamma$, but this does not occur.

\begin{figure}[h!]
\centerline{\hbox{\psfig{figure=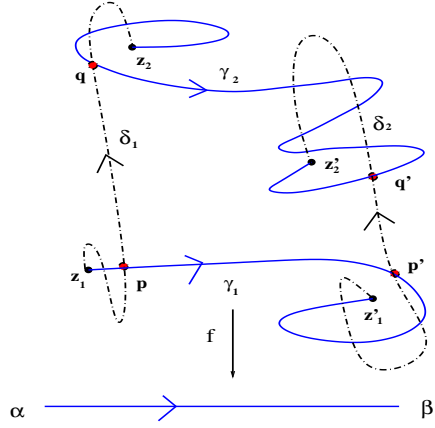,height=2.2in,width=2.2in}}}
\caption{The  points $p$, $p'$, $q$, $q'$ and the curves $\delta_1$, $\delta_2$, $\gamma_1$ and $\gamma_1$}
\label{fig1}
\end{figure}

Take $\gamma_1$ and $\gamma_2$. Since $G_1$ is a domain we can join 
$z_1$ to $z_2$ by an arc $\delta_1$ in $G_1$ and similarly $z_1'$ to $z_2'$ 
by an arc $\delta_{2}$ in $G_2$. 
If $\delta_2$ is oriented from $z_1'$ to $z_2'$, let
$p'$ be the point where, for the last  time, $\gamma_1$ meets $\delta_2$
 and $q'$ be the point where, for the first time, $\gamma_2$ meets $\delta_2$.
If $\delta_1$  is oriented from $z_1$ to $z_2$, let 
$p$ be the point where, for the last time, $\gamma_1$ 
meets $\delta_1$ and $q$ be the point where, for the first  time, 
$\gamma_2$ meets $\delta_1$, these might look like Figure \ref{fig1}.

\begin{figure}[h!]
\centerline{\hbox{\psfig{figure=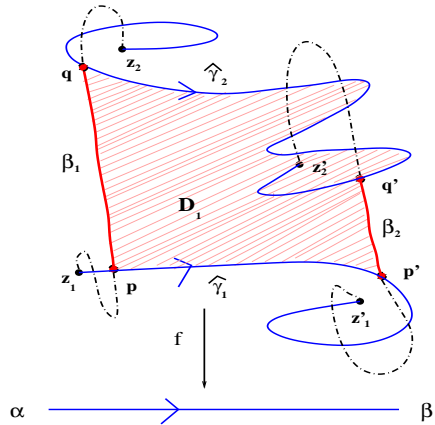,height=2.2in,width=2.2in}}}
\caption{The arcs $\beta_1$, $\beta_2$, $\widehat{\gamma_1}$ , $\widehat{\gamma_2}$  and $D_1$}
\label{fig2}
\end{figure}

Now  we denote by $\beta_1$ the part of $\delta_1$  which joins 
the points $p$ and $q$, by  
$\beta_2$ the part of $\delta_2$ which 
joins the points $p'$ and $q'$, by $\widehat{\gamma_1}$ the part of 
$\gamma_1$ which joins the points $p$ and $p'$,  oriented from $p$ to $p'$, 
and by $\widehat{\gamma_2}$ the part of $\gamma_2$  which joins 
the points  $q$ and $q'$, oriented from $q$ to  $q'$. Then
$\widehat{\gamma_1} \beta_2 \widehat{\gamma_2}^{-1} \beta_{1}^{-1}$  
is a simple closed  curve with an interior $D_1$, see Figure \ref{fig2}.  \\

\noindent {\bf Step 2. The map on the Disc $D_{1}$.} 

Recall that the disc $D_{1}$ has boundary ${\beta_{1}} \cup {\beta_{2}}\cup {\widehat{\gamma_1}} \cup {\widehat{\gamma_2}}$, the end points of the curve  ${\beta_{1}}$  are the points  $p$ and $q$, the end points of  the curve ${\beta_{2}}$ are the points  $p'$ and $q'$, the end points of the curve ${\widehat{\gamma_1}}$ are the points $p$ and $p'$ and the end points  of  the curve ${\widehat{\gamma_2}}$ are the points $q$ and $q'$,  see  Figure \ref{fig2}. The function $f$ maps ${\widehat{\gamma_i}}$ inyectively into the cut $\gamma$,  for $i=1,2$, and we consider $f({\beta_{1}})$ and $f({\beta_{2}})$  two non intersecting curves (with possible self intersections)  with   ends at $f(p),f(q)$  and $f(p'),f(q')$  respectively.

A natural question arises: Where is mapped $D_1$ under $f$?

Observe that $f(D_1) $ can be  either unbounded or bounded. If $f(D_{1})$ is unbounded,  so there is a pole in $D_{1}$. Thus we ruled out this case. Necessarily $f(D_{1})$ must be   bounded and $f(\beta_{1})$ or $f(\beta_{2})$ need not be closed curves. This is the  missing  case  in Baker's proof. \\

\begin{figure}[h!]
\centerline{\hbox{\psfig{figure=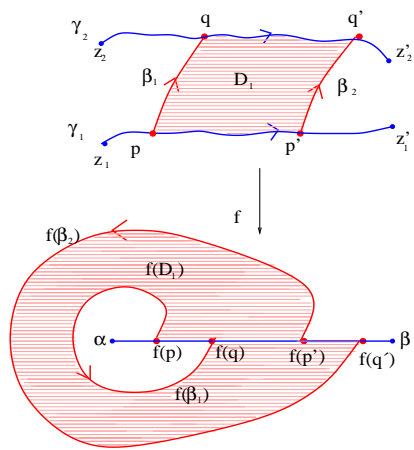,height=2.3in,width=2.1in}}}
\caption{ (a) $\gamma_1$ and $\gamma_2$ have the same orientations}
\label{fi8}
\end{figure}

\begin{figure}[h!]
\centerline{\hbox{\psfig{figure=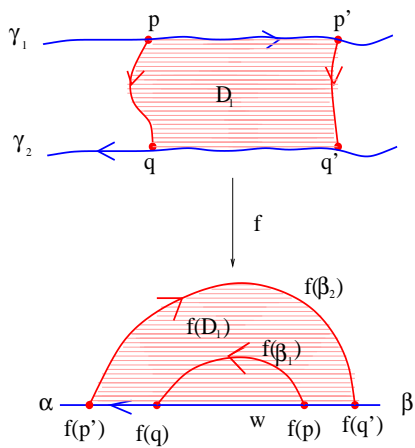,height=2.3in,width=2.1in}}}
\caption{ (b) $\gamma_1$ and $\gamma_2$  have opposite orientations }
\label{fi9}
\end{figure}

Remember that, the orientations of $\gamma_{1}$ and $\gamma_{2}$ are given by the chosen orientation in $\gamma$ as in Step 1 above.  Two main possibilities arises when we consider the orientation of $\gamma_{1}$  together  with the order of the  set of points $\{p, p'\}$  and the orientation of  $\gamma_{2}$  together with  the order of the  set  of points $\{q, q'\}$. Let  us define $a<b$ for $a, b$ points in an oriented curve  $\gamma (t)$, if $\gamma(t_{1})=a$ and $\gamma(t_{2})=b$ and $t_{1}<t_{2}$.  The possibilities are: (a) $\gamma_{1}$ and $\gamma_{2}$ preserve the same order, that is, if $p<p'$,  then $q<q'$,  see Figure \ref{fi8}, or (b) $\gamma_{1}$ and $\gamma_{2}$ reverse the order, that is, $p < p'$ but $q > q'$ or $q < q'$ but $p > p'$, see Figure \ref{fi9}.

\begin{figure}[h!]
\centerline{\hbox{\psfig{figure=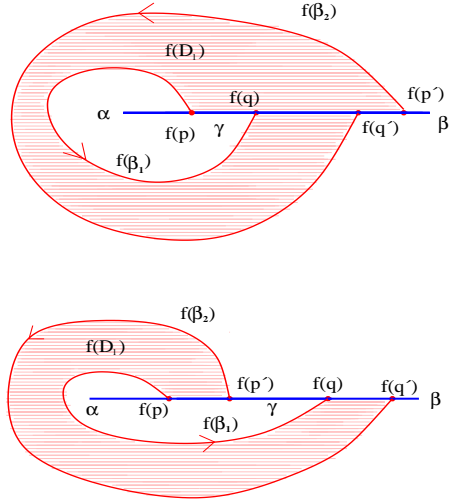,height=2.6in,width=2.3in}}}
\caption{$f({\widehat{\gamma_1}} )\cap f({\widehat{\gamma_2}}) \neq \emptyset$ and  $f({\widehat{\gamma_1}} )\cap f({\widehat{\gamma_2}}) =  \emptyset$}
\label{fig4}
\end{figure}

On the other hand, the curves $f({\beta_{i}})$ has winding number  either $+1$, $0$ or $-1$ with respect to the points $\alpha$ and $\beta$.  So, several possibilities occurs for the topology of $f(D_{1})$  accordingly to how the intervals $f( {\widehat{\gamma_1}})$ and $f( {\widehat{\gamma_2}})$ are placed in the cut $\gamma$.   In Figure \ref{fig4}, there are two examples, one when $f( {\widehat{\gamma_1}}) \cap f( {\widehat{\gamma_2}}) \neq \emptyset$ and the other when $f( {\widehat{\gamma_1}}) \cap f( {\widehat{\gamma_2}})  = \emptyset$. \\

In the next step we  will finish the proof of the theorem. The construction indicated  will work as well for    (a) and (b), but is important to know that the difference exist.\\

\noindent {\bf Step 3. Unbounded curves and their pre-images.}

From now  on, we will assume without lost of generality that $f({\beta}_{2})$ surrounds $f({\beta}_{1})$,  it may look like Figures \ref{fi9} or \ref{fig4}.  Also  we assume that $\gamma_{1}$, $\gamma_{2}$ and $\gamma$ are compatibly oriented, as in Step 1. 

For $x,w \in {\C}$, we  denote by $\overline{xw}$ the oriented segment from $x$ to $w$ and  by $T_{x}(\tau)$  the tangent at $x$ of some parametrization of a curve $\tau$.

\begin{figure}[h!]
\centerline{\hbox{\psfig{figure=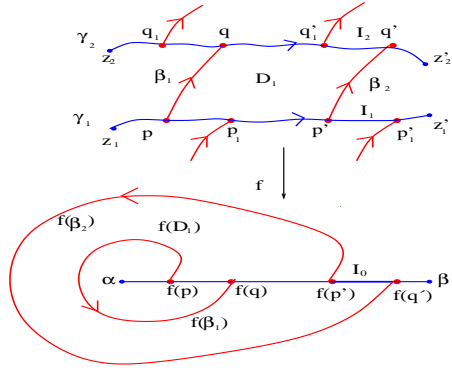,height=1.9in,width=2.3in}}}
\caption{ The pre-images of $f(p), f(q) f(p')$ and $f(q')$}
\label{fig8}
\end{figure}

The step consists of considering certain unbounded curves on the regions $G_{1}$ and $G_{2}$ and their pre-images.  We recall that ${\beta}_{1} \subset G_{1}$ and $\beta_{2} \subset G_{2}$. The curve $f(\beta_{1})$ has end points at $f(p)$ and $f(q)$ in $\gamma$, and the curve $f(\beta_{2})$ has end points at $f(p')$ and $f(q')$ in $\gamma$. So we define their pre-images on $\gamma_{1}$ and on $\gamma_{2}$ as follows: $f^{-1}(f(q')) \cap {\gamma}_{1}=p_{1}'$, $f^{-1}(f(q')) \cap {\gamma}_{2}=q'$, $f^{-1}(f(p')) \cap {\gamma}_{1}=p'$, $f^{-1}(f(p')) \cap {\gamma}_{2}=q_{1}'$, $f^{-1}(f(p)) \cap {\gamma}_{1}=p$,  $f^{-1}(f(p)) \cap {\gamma}_{2}=q_{1}$ and $f^{-1}(f(q)) \cap {\gamma}_{1}=p_{1}$, $f^{-1}(f(q)) \cap {\gamma}_{2}=q$, see Figure \ref{fig8} as an example. The point $p_{1}'$ is the end point of some pre-image of $f(\beta_{1})$, $q_{1}'$ is the beginning of another pre-image of $f(\beta_{1})$, $p_{1}$ is the end point of some pre-image of $f(\beta_{2})$ and $q_{1}$ is the beginning of another pre-image of $f(\beta_{2})$.

For brevity,  we  define $I_{1}$ as the interval $\overline{p'p'_{1}}$ and $I_{2}$  as the interval $\overline{q'_{1}q'}$. Thus $I_{1}$ and $I_{2}$ are pre-images of the interval $I_{0}=\overline{f(p')f(q')}$ in $\gamma_{1}$ and $\gamma_{2}$ respectively. We have   two situations.

\begin{figure}[h!]
\centerline{\hbox{\psfig{figure=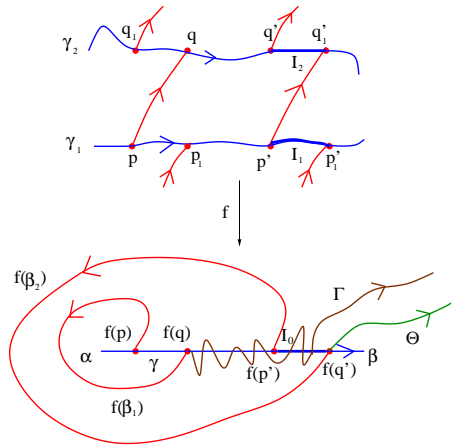,height=2.3in,width=2.3in}}}
\caption{The curves $\Gamma$ and $\Theta$}
\label{fig100}
\end{figure}

(i) Let us consider an unbounded oriented curve ${\Gamma} \subset G_{1}$ beginning at $f(q)$, and an unbounded oriented curve ${\Theta} \subset G_{2}$ beginning at $f(q')$, as  in  the   Observation  1 in Section 1, see for instance  Figure \ref{fig100}.  More conveniently, we are interested in the piece of such curves  complementary to $f(D_{1})$. Denote by $B$ the complement of $f(D_{1})$ in the sphere and let $\Gamma_{0}={\Gamma} \bigcap B$ and  $\Theta_{0}={\Theta} \bigcap B$.

\begin{figure}[h!]
\centerline{\hbox{\psfig{figure=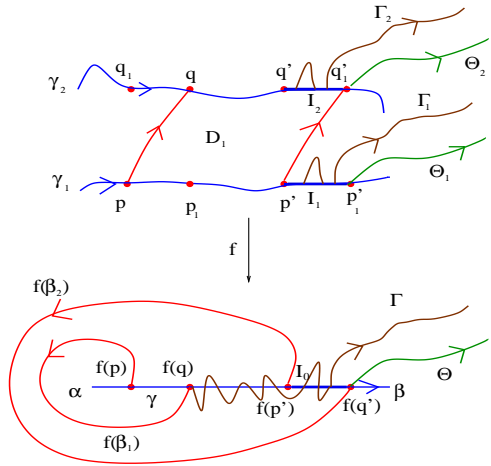,height=2.4in,width=2.5in}}}
\caption{ The  curves  $\Gamma_{1}$, $\Theta_{1}$, $\Gamma_{2}$ and $\Theta_{2}$ }
\label{fig99}
\end{figure}

(ii) The curves $\Gamma$ and $\Theta$ may oscillate and may intersect the interval $I_{0}$ in many  points, so in this case $\Gamma_{0}$ and $\Theta_{0}$ are a union of curves beginning at points in $I_{0}$. By applying the Kaplan's theorem to each of these curves,  we  consider their pre-images beginning at points in $I_{1}$, denoted by $\Gamma_{1}$ and $\Theta_{1}$ respectively and pre-images beginning at $I_{2}$, denoted by $\Gamma_{2}$ and $\Theta_{2}$ respectively, none of these curves intersects $D_{1}$ or more generally $f^{-1}f(D_{1})$, it may look like Figure  \ref{fig99}. If $N_{\infty}$ is any neighborhood of infinity, we have  $\Gamma_{i} \bigcap N_{\infty} \neq \emptyset$   also  $\Theta_{i} \bigcap N_{\infty} \neq \emptyset$, $i=0,1,2$, they are unbounded.

We have now  two cases,  either ( A ) the intersection of the set $\Gamma_{1}$  or $\Theta_{1}$  with $I_{1}$ is finite, consequently the same for $\Gamma_{2}$   or $\Theta_{2}$, or  (B) the intersection of both sets is not finite.

\begin{figure}[h!]
\centerline{\hbox{\psfig{figure=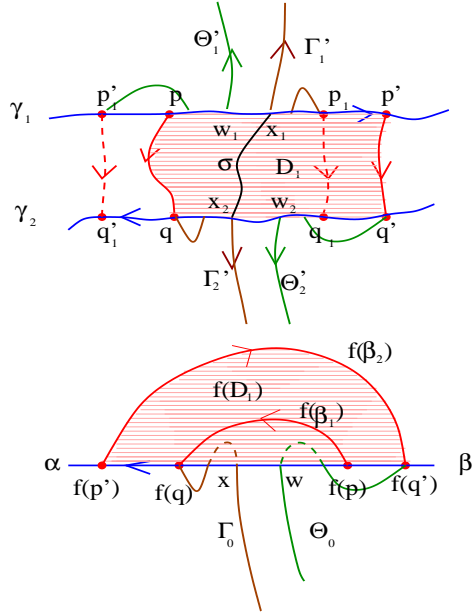,height=3.2in,width=2.4in}}}
\caption{ The  curves  $\Gamma_{0}$, $\Theta_{0}$ and the points $z_{i}$ and $w_{j}$, $i=0,1,2$.}
\label{fig999}
\end{figure}

Case A. Assume without lost of generality that $\Gamma_{i}$ intersects $I_{i}$ in a finite set, $i=1,2$. We consider the component of $\Gamma_{1}$ which is unbounded and denote it by $\Gamma_{1}'$, similarly we have an unbounded component $\Gamma_{2}'$. Both curves are in $G_{1}$ and recall that their closure in the sphere is a continua that contains infinity. Let us denote  $x_{1}=\Gamma_{1}' \bigcap I_{1}$ and $x_{2}=\Gamma_{2}' \bigcap I_{2}$, observe that  the pairs  $(T_{x_{1}}(I_{1}),T_{x_{1}}(\Gamma'_{1}))$ and $(T_{x_{2}}(I_{2}),T_{x_{2}}(\Gamma'_{2}))$ are sent conformally by $f'$ to the corresponding pair $(T_{f(x_{1})}(I_{0}), T_{f(x_{1})}(\Gamma_{0}))$. Under such conditions, for any path $\sigma$  joining $x_{1}$ with $x_{2}$ which does not intersect $\Theta_{1}$ nor $\Theta_{2}$,  then the curve $-\Gamma_{1}' \bigcup \sigma \bigcup \Gamma_{2}$ disconnects $\Theta_{1}$ from $\Theta_{2}$, it may look like Figure \ref{fig999}. Therefore, if $\sigma \in G_{1}$, then $G_{2}$ is disconnected.\\

Case B.  In this case the intersection $I_{i} \bigcap \Gamma_{i}$ is an infinite collection of points $\{x^{j}_{i}\}$, equally $I_{i} \bigcap \Theta_{i}$ is an infinite  collection of points $\{w^{j}_{i}\}$ for $i=1,2$ and $j \in \N$. Being $I_{1}$ and $I_{2}$ compact, the sequence $\{x^{j}_{i}\}$ has at least an accumulation point, say $x_{i}$,  and let $w_{i}$ be an accumulation point for the sequence $\{w^{j}_{i}\}$, $i=1,2$. Again we have two situations, either (1) at least one of the points $x_{1}$ or $w_{1}$ is in the Fatou ser, or (2) both points $x_{1}$ and $w_{1}$ are in the Julia set.\\

(1) Assume  without lost of generality that $x_{1}$ and  so $x_{2}$ are in the Fatou set. Consider the closure on the sphere $\overline{\Gamma_{i}}$ of $\Gamma_{i}$, $i=0,1,2$. Let $\sigma$ be a path between $x_{1}$ and $x_{2}$ that does not intersect  $\Theta_{i}$, $i=1,2$.  As in Case A, the pairs $(T_{x_{1}}(I_{1}),T_{x_{1}}(\overline{\Gamma_{1}}))$ and $(T_{x_{2}}(I_{2}),T_{x_{2}}(\overline{\Gamma_{2}}))$ are sent conformally by $f'$ to the corresponding pair $(T_{f(x_{1})}(I_{0}), T_{f(x_{1})}(\overline{\Gamma_{0}}))$, where  $(T_{x_{i}}(I_{i}),T_{x_{i}}(\overline{\Gamma_{i}}))$ means $lim_{x_{i}^{j} \rightarrow x_{i}}(T_{x_{i}}(I_{i}),T_{x_{i}}(\Gamma_{i}))$, $i=0,1,2$.

Under such conditions, for any path $\sigma$  joining $x_{1}$ with $x_{2}$ which does not intersect  neither $\Theta_{1}$ nor $\Theta_{2}$,  the set $\overline{\Gamma_{1}} \bigcup \sigma \bigcup \overline{\Gamma_{2}}$ disconnects $\Theta_{1}$ from $\Theta_{2}$. Therefore, if $\sigma \in G_{1}$, then $G_{2}$ is disconnected.\\

(2) Assume that $x_{i}$ and $w_{i}$ are in the Julia set, $i=1,2$. Consider paths $\sigma_{j}$ between $x_{1}^{j}$ and $x_{2}^{j}$, $j \in {\N}$ and let $\Sigma=\overline{\bigcup_{j}{\sigma}_{j}}$ be the closure of the union of all the paths $\{\sigma_{j}\}$ in the sphere. Also, the pairs $(T_{x_{1}}(I_{1}),T_{x_{1}}(\overline{\Gamma_{1}}))$ and $(T_{x_{2}}(I_{2}),T_{x_{2}}(\overline{\Gamma_{2}}))$ are sent conformally by $f'$ to the corresponding pair $(T_{f(x_{1})}(I_{0}), T_{f(x_{1})}(\overline{\Gamma_{0}}))$. As explained in (ii), the intersection of the sets $\overline{\Gamma_{1}},\overline{\Gamma_{2}},\overline{\Theta_{1}}$ and $\overline{\Theta_{2}}$ with $D_{1}$ is empty.

Observe that the set $\Sigma \bigcup \overline{\Gamma_{1}} \bigcup \overline{\Gamma_{2}}$ disconnects $G_{2}$, since in any neighborhood of $x_{1}$ and $x_{2}$ there are points that belong to $G_{2}$. Now, if all $\sigma_{j} \in G_{1}$, then $\Sigma \bigcap G_{2}=\emptyset$ and $\Sigma \bigcup \overline{\Gamma_{1}} \bigcup \overline{\Gamma_{2}}$ is disjoint of $G_{2}$, therefore in this case $G_{2}$ is disconnected.\\

 In all these cases $G_{2}$ is disconnected which is a contradiction.  Thus we have finished the proof of  the Theorem \ref{ba}. \\

\noindent {\bf Remark.} This  above proof applies also to the case of a transcendental meromorphic map with a finite number of poles, see \cite{pat}, since we can choose a disc $D_{1}$  without poles exactly as in the proof  of the Theorem and this  case proceeds as above.\\

\noindent {\bf Acknowledgmets.} The  authors would like to thank to P. Rippon, G. Stallard,  M. Montes de Oca and the members of the holomorphic  dynamics seminar in UNAM,  for  their comments and support when we were having different drafts  and ideas of the  the proof. We thank  specially to J.  Kotus  for her patience  to  listening  our arguments and  also for very interesting discussions of the proof.


\end{document}